\documentclass[twoside,a4paper,reqno,11pt]{amsart} 
\usepackage{amsfonts, amsbsy, amsmath, amssymb, latexsym}
\usepackage{mathrsfs,array}
\usepackage[top=30mm,right=30mm,bottom=30mm,left=30mm]{geometry}
\usepackage{stmaryrd}
\usepackage{bm}

\usepackage{hyperref}

            \def\e{\epsilon}

            \def\a{\alpha} 
             \def\b{\beta} 
            \def\c{\chi} 
             
            \def\e{\epsilon} 
            \def\O{\Omega} 
            \def\o{\omega} 
            \def\l{\lambda}

            \def\pf{\noindent {\bf Proof } }

            \def\F{\mathbb{F}} 
    \def\W{\mathcal{W}}

            \def\Z{\mathbb{Z}}

            \def\d{\delta}

     \def\G{\Gamma} 
             
\def\no{\noindent}

            \newtheorem{thm}{Theorem}[section] 
            \newtheorem{prop}[thm]{Proposition} 
            \newtheorem{lem}[thm]{Lemma} 
             
             \newtheorem{cor}[thm]{Corollary}

\begin{document}
\author{Daniele Dona}
\address{Alfr\'ed R\'enyi Institute of Mathematics, Re\'altanoda utca 13-15, H-1053, Budapest, Hungary}
\email{dona@renyi.hu}

\author{Martin W. Liebeck} 
\address{Imperial College London, 180 Queen's Gate, South Kensington, London, SW7 2AZ, United Kingdom}
\email{m.liebeck@imperial.ac.uk}

\author{Kamilla Rekv\'enyi} 
\address{Imperial College London, 180 Queen's Gate, South Kensington, London, SW7 2AZ, United Kingdom}
\email{k.rekvenyi19@imperial.ac.uk}

           \title{Involutions in finite simple groups as products of conjugates}

           \begin{abstract}

Let $G$ be a finite non-abelian simple group, $C$ a non-identity conjugacy class of $G$, and $\G_C$ the Cayley graph of $G$ based on $C \cup C^{-1}$. Our main result shows that in any such graph, there is an involution at bounded distance from the identity.
               
           \end{abstract}

\date{\today}
\maketitle

\section{Introduction} \label{intr}

For a finite group $G$ with a generating set $S$ that is closed under inversion, the Cayley graph $\G(G,S)$ is defined to be the graph with vertex set $G$ and edge set consisting of all $\{g,gs\}$ for $g\in G,s\in S$. This paper is a contribution to the study of the Cayley graphs of finite (non-abelian) simple groups $G$ in the case where the generating set is of the form $C \cup C^{-1}$, where $C$ is a non-identity conjugacy class of $G$. By \cite{LSh1}, the diameter of such a Cayley graph is bounded above by $c\log |G|/\log |C|$, where $c$ is an absolute constant. Also, for $G$ a simple group of Lie type of rank $l$, an explicit upper bound, linear in $l$, for the diameters of all such Cayley graphs, is obtained in \cite{LL}. 

Here we study a more refined question concerning distances in these Cayley graphs $\G_C:= \G(G,C\cup C^{-1})$ for simple groups $G$. For a subset $H \subseteq G$ not containing the identity, define a number $d(H)$ as follows. Let $d_C$ denote the distance function on $\G_C$, and set
\[
d_C(H) = \hbox{min}\left(d_C(1,h) : h \in H\right),
\]
so that $d_C(H) = \hbox{min}\left( k : (C\cup C^{-1})^k \cap H \ne \emptyset \right)$. Now define
\[
d(H) = \hbox{max}_{C\neq 1} \;d_C(H).
\]
In other words, for every non-identity conjugacy class $C$ of $G$, we have $(C\cup C^{-1})^k \cap H \ne \emptyset$ for some $k \le d(H)$.

Clearly, there are some subsets $H$ for which $d(H)$ is equal to $d_G : = \hbox{max}_{C\neq 1} \,\hbox{diam}(\G_C)$, hence grows with the rank of $G$. Our main results show that this is not the case when $H= Inv(G)$, the set of involutions in the simple group $G$. Indeed, we show that $d(Inv(G))$ is bounded above by an absolute constant: explicitly, by 2 for alternating groups $G = A_n\,(n\ge 6)$, by 3 for $G=A_{5}$, and by 12 for classical groups. Constants for $G$ an exceptional group of Lie type and $G$ sporadic already follow from existing results: see Proposition~\ref{excspor}. Thus for any non-identity conjugacy class $C$ of a finite simple group $G$, there is a bounded product of elements of $C\cup C^{-1}$ that is an involution.

%In this paper we consider a refinement of width questions in finite simple groups. 

%Let $T$ be a non-abelian simple group, and $S$ a generating set of $T.$ Define the width of $T$ with respect to $S,$ denoted $w_S(T),$ to be the minimal $k\in \mathbb{N}$ such that any element of $T$ can be expressed as a product of at most $k$ elements of $S.$  Finding this number has been considered with respect to several generating sets, for example the case where $S$ is a conjugacy class, also referred to as conjugacy width, the case where $S$ is the union of a conjugacy class and its inverse, also referred to as inverse conjugacy width, denoted $c_i(T)$ or the case when $S$ is the set of involutions in $T.$ 
%It was shown in \cite{alexmalcolm} that we can express every element in a finite simple group as a product of at most $4$ involutions. 

%Here we consider an opposite problem, and provide a number $l$ such that for any non-central element $g$ in $T$ we can express an involution as a product of $l$ conjugates of $g$ or $g^{-1}.$  We will call this number $i(T). $ Note that this is a lower bound for %$c_i(T).$

We now state our results. For classical groups, our proofs are inductive, which necessitates consideration of the general linear, unitary and orthogonal groups rather than just the simple groups. For $\e = \pm$, $GL_n^\e(q)$ denotes $GL_n(q)$ if $\e=+1$, and $GU_n(q)$ if $\e = -1$, with similar notation for $SL_n^\e(q)$. Also when $n$ is even, $O_n^\e(q)$ denotes the orthogonal group $O_n^\pm(q)$ as usual, and when $n$ is odd it denotes the orthogonal group $O_n(q)$.

\begin{thm}\label{glgu} Let $G = GL_n^\e(q)$ and $S = SL_n^\e(q)$, with $n\ge 2$ and $(n,q,\e) \ne (2,2,\e)$ or $(3,2,-)$. Let $g \in G\setminus Z(G)$, and $C  = g^S$. Then there exists an involution $t \in S$, with $t \not \in Z(S)$ when $n \ne 2$, and elements $x_i \in C$, $y_i \in C^{-1}$ for $1\le i \le 6$, such that $\prod_{i=1}^6 x_iy_i = t$.
\end{thm}

\begin{thm}\label{sp} Let $G = Sp_{2m}(q)$, with $m\ge 2$. Let $g \in G\setminus Z(G)$, and $C  = g^G$. Then there exists an involution $t \in G\setminus Z(G)$, and elements $x_i \in C$, $y_i \in C^{-1}$ for $1\le i \le 6$, such that 
$\prod_{i=1}^6 x_iy_i = t$.
\end{thm}

\begin{thm}\label{othm} Let $G = O_n^\e (q)$ and $S = \O_n^\e(q)$, with $n\ge 7$, and $q$ odd if $n$ is odd. Let $g \in G\setminus Z(G)$, and $C  = g^S$. Then there exists an involution $t \in S\setminus Z(S)$, and elements $x_i \in C$, $y_i \in C^{-1}$ for $1\le i \le 6$, such that $\prod_{i=1}^6 x_iy_i = t$.
\end{thm}

The following corollary follows directly from these theorems, except in the case where $G = PSL_2(q)$ with $q$ odd, in which case $C^3=G$ for all non-identity conjugacy classes $C$, by \cite[Thm 4.2 (a)]{arad1}. 

\begin{cor}\label{dcor} Let $G$ be a classical simple group, and let $1 \ne g \in G$. Then there are $12$ conjugates of $g^{\pm 1}$ with product equal to an involution. Thus $d(Inv(G)) \le 12$.
\end{cor} 

In Section \ref{alg} we provide a different proof that for $G=PSL(n,q)$ and $1 \ne g \in G$, there is a product of at most $48$ conjugates of $g^{\pm 1}$ that is equal to an involution (see Theorem \ref{sl}). Although this bound 48 is worse than the bound in Corollary \ref{dcor}, the proof is constructive, in the sense that it provides an algorithm that, given $g$, provides an explicit product of conjugates of $g$ that is an involution.

%\no {\it Remark } For exceptional groups of Lie type, there is an upper bound of $8(5r+7)$ given by \cite{LL}, where $r$ is the Lie rank and for sporadic groups there is an upper bound of $6,$ as the covering numbers for those are at most $6.$

For alternating groups we prove the following.

\begin{thm}\label{alt}
     For all $g\in S_n \setminus 1$, there exists $h\in A_n$ such that $[g,h]$ has order 2, unless $n=5$ and $g$ is a $5$-cycle. Thus $d(Inv(A_n)) = 2$ for $n\ge 6$, and $d(Inv(A_5)) = 3$.
\end{thm}

%Note that since $cn(A_5)=3$ we can express an involution as a product of $3$ conjugates of a $5$-cycle in $S_5.$

For the remaining simple groups, we note the following result, which follows directly from \cite[Thm. 2]{LL} (for exceptional groups of Lie type, including the Tits group $^{2}F_{4}(2)'$) and \cite{Zisser} (for sporadic groups).

\begin{prop}\label{excspor} If $G$ is a simple group of exceptional Lie type, then $d(Inv(G)) \le 376$; and if $G$ is a sporadic simple group, then $d(Inv(G)) \le 6$.
\end{prop}

Finally, we mention a consequence of our bounds on $d(Inv(G))$ for simple groups $G$. Namely, we obtain a new bound on the orbital diameters of primitive permutation groups of simple diagonal type, improving a result in \cite{kamilla}. See Section \ref{orbdiam} for details.

The layout of the paper is as follows. In Sections \ref{secglgu}--\ref{secothm} we prove Theorems \ref{glgu}--\ref{othm}, and Section \ref{alg} contains our constructive proof for $PSL_n(q)$. The proof of Theorem \ref{alt} can be found in Section \ref{secalt}, and the final section contains our result on orbital diameters.

%As a corollary to our result, we can improve the bound in \cite{kamilla} concerning the orbital diameter of primitive groups of simple diagonal type. 

%\begin{coroll}
%Let $k\geq 3$ and let $T$ be a finite non-abelian simple group. Then $$orbdiam(T^k.S_k)\leq \min \begin{cases}
 %   24(k-1)i(T)c_i(T)\\
 %   24(k-1)c_i^2(T)
%\end{cases}.$$

%\end{coroll}

\section{Proof of Theorem \ref{glgu}} \label{secglgu}

The proof will be inductive, and the following lemma is the $n=2$ case, which forms the base for the induction.

\begin{lem}\label{base} Let $G = GL_2^\e(q)$ and $S = SL_2(q)$, with $q> 2$. If $g \in G\setminus Z(G)$, and $C = g^S$, then 
there exist elements $x_i \in C$, $y_i \in C^{-1}$ for $1\le i \le 6$, such that 
\[
\prod_{i=1}^6 x_iy_i  = \left\{ \begin{array}{l} -I,\hbox{ if }q \hbox{ is odd} \\ \begin{pmatrix} 1&1 \\ 0&1 \end{pmatrix}, \hbox{ if }q \hbox{ is even} \end{array} \right.
\]
\end{lem}

\pf First choose $h \in S$ such that $[g,h] \ne I$, and note that $[g,h] \in CC^{-1} \subseteq S$.
If $q$ is even, then $(CC^{-1})^3 = S$ by \cite[Thm 4.2 (a)]{arad1}, and the conclusion follows.

Now assume that $q$ is odd, and let $x = [g,h]$. If $x = -I$, there is nothing to prove, so suppose $x \ne -I$. By \cite[p.43]{arad1}, the number of solutions to the equation $x_1\cdots x_6 = -I$ with each $x_i \in x^S$ is equal to
\[
\frac{|x^S|^6}{|S|} \sum_{\chi \in Irr(S)} \frac{\chi(x)^6\chi(-I)}{\chi(1)^{5}}.
\]
Using the character table of $S$ (see for example \cite[p.229]{dorn}), we check that this sum is positive, and the result follows. 
$\;\;\;\Box$

\vspace{4mm}
\no {\bf Proof of Theorem \ref{glgu}}

\vspace{2mm} We prove the theorem first for $G = GL_n(q)$, and then indicate the changes needed to handle $GU_n(q)$. 
Let $G = GL_n(q)= GL(V)$ and $S = SL_n(q)$, with $(n,q) \ne (2,2)$, let $g \in G\setminus Z(G)$, and $C  = g^S$. The proof goes by induction on $n$. The base case $n=2$ is given by Lemma \ref{base}, so assume $n\ge 3$. Assume also that 
$(n,q)$ is not $(3,2)$, $(3,3),$ $(3,4)$, $(3,5)$ $(4,2)$ or $(4,3)$. In these cases the result follows by calculations in GAP \cite{gap}.

\vspace{2mm}
\no {\bf (1) } Consider first the case where $g$ is {\it decomposable}, i.e. $g$ lies in a proper subgroup $GL_{n_1}(q) \times GL_{n_2}(q)$ of $G$ preserving a decomposition $V = V_1 \oplus V_2$, where $\dim V_i = n_i$. Say $g = g_1\oplus g_2$, where $g_i \in GL(V_i)$. 

Assume that there exists $i$ such that $g_i \not \in Z(GL(V_i))$, say $i=1$. Then by induction, one of the following holds:
\begin{itemize}
\item[(i)] there is an involution $t_1 \in SL(V_1)$, and $SL(V_1)$-conjugates $x_j'$ and $y_j'$ ($1\le j\le 6$) of $g_1$ and $g_1^{-1}$ respectively, such that $\prod_{1}^6 x_j' y_j' = t_1$;
\item[(ii)] $(n_1,q) = (2,2)$.
\end{itemize}
In case (i), if we let $x_j = x_j' \oplus g_2$ and $y_j = y_j'\oplus g_2^{-1}$, then these are $S$-conjugates of $g$ and $g^{-1}$, with product $t_1 \oplus 1$, a non-central involution in $S$. In case (ii) we have $n\ge 5$ by assumption, so $n_2\ge 3$. If $g_2 \ne 1$, work as above with $g_2$ instead of $g_1$; and if $g_2 = 1$, we can replace $g_1$ by $(g_1,1) \in GL_3(2)$ and proceed as above.

Now assume that $g_i \in Z(GL(V_i))$ for $i=1,2$, so $g = \l I_{n_1} \oplus \mu I_{n_2}$ for some $\l,\mu \in \F_q^*$. Then we can replace $g$ by the $S$-conjugate $h_1\oplus h_2$, where $h_1 = (\l I_{n_1-1},\mu)$, $h_2 = (\l,\mu I_{n_2-1})$, and proceed as in the previous paragraph.

\vspace{2mm}
\no {\bf (2) } Now consider the case where $g$ is indecomposable, i.e. does not preserve any proper decomposition $V = V_1 \oplus V_2$. Let $g=su$ be the Jordan decomposition of $g$, where $s$ is semisimple, $u$ is unipotent, and $us=su$. The structure of the centralizer $C_G(s)$ is well known, and we have 
\[
u \in C_G(s) = \prod_{i=1}^k GL_{m_i}(q^{d_i}) \le \prod_{i=1}^k GL_{m_id_i}(q) \le G,
\]
where $\sum m_id_i = n$. Since $g$ is indecomposable, we have $k=1$, so $C_G(s) = GL_m(q^d)$ with $md=n$, and also $u$ is a single Jordan block $J_m \in GL_m(q^d)$. Moreover, if $2<m<n$, or if $m=2$ with $q$ even, then we can apply induction in the group $GL_m(q^d)$ to obtain the conclusion of the theorem, so assume that $m=1,\,2$ or $n$ (with $q$ odd if $m=2$).

If $m=1$ then $g=s$ and $C_G(g) = GL_1(q^n)$; if $m=n$ then $g = \l J$ with $J$ a unipotent Jordan block, and $C_G(g)$ has order $(q-1)q^{n-1}$; and if $m=2$ then $C_G(g) = C_{GL_2(q^d)}(J_2)$ has order $(q^{d}-1)q^d$. In all cases,
\begin{equation}\label{cbd}
|C_G(g)| < q^n.
\end{equation}

Suppose first that $m=n$, so $g = \l J$. Let $C = J^S$, the conjugacy class of $J$ in $S = SL_n(q)$, and let $t$ be a non-central involution in $S$. By \cite[p.43]{arad1}, the number of solutions to the equation $\prod_1^6x_iy_i = t$ with $x_i \in C$, $y_i \in C^{-1}$ is
\begin{equation}\label{csum}
N:= \frac{|C|^{12}}{|S|} \sum_{\c \in Irr(S)} \frac{\c(J)^6\c(J^{-1})^6\c(t)}{\c(1)^{11}}.
\end{equation}
To show that this sum is nonzero, we use the following facts; in (c), $k(S)$ denotes the number of conjugacy classes in $S$:
\begin{itemize}
\item[(a)] $|\c(J)| = |\c(J^{-1})| \le q^{n/2}$ for all $\c \in \hbox{Irr}(S)\setminus 1$, by (\ref{cbd});
\item[(b)] $\c(1) \ge \frac{q^n-1}{q-1}-1 > q^{n-1}$ for all $\c \in \hbox{Irr}(S)\setminus 1$, by \cite[Thm. 3.1]{TZ};
\item[(c)] $k(S) \le 2.5q^{n-1}$, by \cite[Prop. 3.6]{FG}.
\end{itemize}
It follows from (a)-(c) that 
\[
\sum_{\c \in Irr(S)\setminus 1} \frac{|\c(J)^6\c(J^{-1})^6|}{\c(1)^{10}} \le \frac{2.5q^{n-1}(q^{n/2})^{12}}{q^{10(n-1)}} = 2.5q^{-3n+9}.
\]
This is less than 1 when $n\ge 4$, so the sum in (\ref{csum}) is nonzero, and the conclusion follows in this case. For $n=3,$ the character table of $SL_3(q)$ is available in \cite{SSF}. We can replace the bound for  $|\c(J)| = |\c(J^{-1})|$ by $2q+1$  for all $\c \in \hbox{Irr}(S)\setminus 1$, so the sum in (\ref{csum}) is nonzero for $q\geq 7$ and the conclusion follows. 

Now suppose that $m=1$, so  $g=s$ and $C_G(g) = GL_1(q^n)$. The determinant map $C_G(g) \mapsto \F_q^*$ is surjective, so we have $g^G = g^S$. Hence, as above, it is sufficient to prove that $N'>0$, where
\begin{equation}\label{ndash}
N':= \sum_{\c \in Irr(G)} \frac{\c(g)^6\c(g^{-1})^6}{\c(1)^{10}}.
\end{equation}
The contribution of the $q-1$ linear characters of $G$ to $N'$ is $q-1$. For the nonlinear characters $\c \in \hbox{Irr}(G)$, again using \cite{FG, TZ}, we have 
\begin{itemize}
\item[(a)] $|\c(g)| = |\c(g^{-1})| \le q^{n/2}$;
\item[(b)] $\c(1) \ge \frac{q^n-1}{q-1}-1 > q^{n-1}$;
\item[(c)] $k(G) \le q^{n}$.
\end{itemize}
Hence the contribution of the nonlinear characters to $N'$ is of absolute value less than 
\[
\frac{q^n(q^{n/2})^{12}}{q^{10(n-1)}} = q^{-3n+10},
\]
which is less than 1 for $n\geq 4$. For $n=3,$ the result follows as again by \cite{Steinberg} we can replace the bound for  $|\c(g)| = |\c(g^{-1})|$ by $6$  for all nonlinear $\c \in \hbox{Irr}(G).$

Finally, suppose $m=2$ with $q$ odd. Here $C_G(g) =  C_{GL_2(q^d)}(J_2)$ and the image of  the determinant map $C_G(g) \mapsto \F_q^*$ has index 2 in $\F_q^*$. It follows that $g^S = g^{G^*}$, where $G^*$ is of index 2 in $G$. The argument of the previous paragraph goes through with $G^*$ replacing $G$, and replacing the inequality in (c) by $k(G^*) \le \frac{1}{2}(q^n+3q^{n-1})$ (see \cite[Cor. 3.7]{FG}).

\vspace{2mm}
This completes the proof of the theorem for $G = GL_n(q)$. 

\vspace{4mm}
We now indicate the changes in the proof needed for $G = GU_n(q)$. First, we need to deal computationally with some small cases, namely 
\begin{equation}\label{exc}
n=3\,(q\le 7), \,n=4\,(q\le 3), \hbox{ and } (n,q) = (5,2),(6,2).
\end{equation}
%({\bf Please check}).
We now assume none of these cases holds. If $g$ is decomposable, i.e. $g$ lies in a proper subgroup $GU_{n_1}(q) \times GU_{n_2}(q)$ of $G$ preserving an orthogonal decomposition $V = V_1 \oplus V_2$, then we argue by induction in similar fashion to the $GL_n(q)$ case, as follows. Write $g = g_1 \oplus g_2$ with $g_i \in GU(V_i)$, and assume first that there exists $i$ such that $g_i \not \in Z(GU(V_i))$, say $i=1$. Then by induction, one of the following holds:
\begin{itemize}
\item[(i)] there is an involution $t_1 \in SU(V_1)$, and $SU(V_1)$-conjugates $x_j'$ and $y_j'$ ($1\le j\le 6$) of $g_1$ and $g_1^{-1}$ respectively, such that $\prod_{1}^6 x_j' y_j' = t_1$;
\item[(ii)] $(n_1,q) = (3,2)$.
\end{itemize}
In case (i), if we let $x_j = x_j' \oplus g_2$ and $y_j = y_j'\oplus g_2^{-1}$, then these are $S$-conjugates of $g$ and $g^{-1}$, with product $t_1 \oplus 1$, a non-central involution in $S$. In case (ii) we have $n\ge 7$ by assumption, so $n_2\ge 4$. If 
$g_2 \not \in Z(GU(V_2))$, work as above with $g_2$ instead of $g_1$; otherwise, $g_2 = \o I_{n_2}$, we can replace $g_1$ by $(g_1,\o) \in GU_4(2)$ and proceed as above.

Now assume that $g_i \in Z(GU(V_i))$ for $i=1,2$, so $g = \l I_{n_1} \oplus \mu I_{n_2}$ for some $\l,\mu \in \F_q^2$. Then we can replace $g$ by the $S$-conjugate $h_1\oplus h_2$, where $h_1 = (\l I_{n_1-1},\mu)$, $h_2 = (\l,\mu I_{n_2-1})$, and proceed as in the previous paragraph.

So assume $g$ is indecomposable, and let $g=su$ be the Jordan decomposition of $g$. Then 
\[
u \in C_G(s) = \prod_{i=1}^k GL_{m_i}(q^{2d_i}) \times  \prod_{i=1}^l GU_{r_i}(q^{e_i}) 
\le \prod_{i=1}^k GL_{2m_id_i}(q) \times \prod_{i=1}^l GU_{r_ie_i}(q)\le G,
\]
where $\sum 2m_id_i + \sum r_ie_i = n$ and all $e_i$ are odd. Since $g$ is indecomposable, $k+l=1$. As in the $GL$ case, we can apply induction to see that one of the following holds:
\begin{itemize}
\item[(i)] $C_G(s) = GU_n(q)$, and $g = \l J$ where $J$ is a single Jordan block in $G$;
\item[(ii)] $C_G(s) = GL_1(q^n)$ with $n$ even, or $GU_1(q^n)$ with $n$ odd;
\item[(iii)] $C_G(s) = GL_2(q^{2d})$ ($n=4d$) or $GU_2(q^e)$ ($n=2e$, $e$ odd), where $q$ is odd, and $u$ is a single Jordan block $J_2 \in C_G(s)$.
\end{itemize}
The maximal possible value of $|C_G(g)|$ is its value in (i), namely $(q+1)q^{n-1}$.

Consider first case (i). We need to argue as in the $GL$ case that the sum in (\ref{csum}) is nonzero. To do this we use the facts:
\begin{itemize}
\item $\c(1) \ge \frac{q^n-q}{q+1}$ for all $\c \in \hbox{Irr}(S)\setminus 1$, by \cite[Thm. 4.1]{TZ};
\item $k(S) \le 8.26q^{n-1}$, by \cite[Prop. 3.10]{FG}.
\end{itemize}
Hence 
\[
\sum_{\c \in Irr(S)\setminus 1} \frac{|\c(J)^6\c(J^{-1})^6|}{\c(1)^{10}} \le \frac{8.26q^{n-1}((q+1)q^{n-1})^{6}}{((q^n-q)/(q+1))^{10}}.
\]
One checks that the right hand side is less than 1 unless $n=3$ or $(n,q)=(7,2),$ $(5,3),$ $(4,5),$ $(4,4)$ or as in the exclusions in (\ref{exc}). For $n=3,$ the character table of $SU_3(q)$ is available in \cite{SSF}. We can replace the bound for  $|\c(g)| = |\c(g^{-1})|$ by $\frac{4q+1}{3}$  for all $\c \in \hbox{Irr}(S)\setminus 1$, so the sum in (\ref{csum}) is nonzero for $q\geq 7$ and the conclusion follows. For $(n,q)=(7,2),$ $(5,3),$ $(4,5),$ $(4,4)$ we can calculate $k(S)$ using GAP, and by replacing $8.26q^{n-1}$ with these exact values we find that the sum in (\ref{csum}) is nonzero. This completes the proof in case (i).

In case (ii) we have $g^S = g^G$ and as before we need to argue that the sum in (\ref{ndash}) is positive. For this we argue as in the previous paragraph, using the bound $k(G) \le 8.26q^n$ from \cite[Prop. 3.9]{FG}, for $n=3$ also the bound $|\c(g)| = |\c(g^{-1})|\le 6$ from \cite{ennola}, and for $(n,q)=(4,4)$ using the fact that $k(G)=470$ from a computation in GAP.  %({\bf Check calcs}) 

Finally, in case (iii) as in the $GL$ case we have $g^S = g^{G^*}$, where $G^*$ is of index 2 in $G$, and we argue in the usual way, using the bound $k(G^*) \le 4.13q^{n-1}(q+1)$ from \cite[Cor. 3.11]{FG} (and also for $(n,q)=(5,2)$ and $(4,4)$, the exact values of $k(G^*)$, computed in GAP). %({\bf Check calcs}) 

This completes the proof of Theorem \ref{glgu}.

\section{Proof of Theorem \ref{sp}}

We prove by induction on $m$ the statement of Theorem \ref{sp}, but including also the case where $m=1$ (which is already proved in Lemma \ref{base}). So let $G = Sp_{2m}(q) = Sp(V)$ with $m\ge 2$, let $g\in G\setminus Z(G)$ and $C = g^G$. Our approach is similar to that in the previous section.

\vspace{2mm}
\no {\bf (1) } Consider first the case where $g$ is decomposable, i.e. $g$ lies in a proper subgroup $Sp_{2m_1}(q) \times Sp_{2m_2}(q)$ of $G$ preserving an orthogonal decomposition $V = V_1 \oplus V_2$, where $\dim V_i = 2m_i$. Say $g = g_1\oplus g_2$, where $g_i \in Sp(V_i)$. Then induction gives the conclusion in exactly the same way as in (1) of the previous section.

\vspace{2mm}
\no {\bf (2) } Now suppose $g$ is indecomposable, and let $g = su$ be the Jordan decomposition. The structure of $C_G(s)$ is as follows: if for $\e \in \{\pm 1\}$, $V_\e$ denotes the $\e$-eigenspace of $s$, then $V_\e$ is non-degenerate and
\[
\begin{array}{ll}
C_G(s) & = Sp(V_1) \times Sp(V_{-1}) \times \prod_{i=1}^k GL_{r_i}(q^{d_i}) \times  \prod_{i=1}^l GU_{s_i}(q^{e_i}) \\
          & \le Sp(V_1) \times Sp(V_{-1}) \times \prod_{i=1}^k Sp_{2r_id_i}(q) \times  \prod_{i=1}^l Sp_{2s_ie_i}(q)\leq G.
\end{array}
\]
Since $g$ is indecomposable, there is only one factor in $C_G(s)$, and so one of the following holds:
\begin{itemize}
\item[(i)] $s = \pm I$, $C_G(s) = G$, and $u$ is an indecomposable unipotent element of $G$;
\item[(ii)] $C_G(s) = GL_r(q^d)$ or $GU_r(q^d)$ with $rd=m$, and $u$ is a single Jordan block $J_r \in C_G(s)$.
\end{itemize}
In case (i), we use \cite[Prop. 2.3, Thm. 3.1]{GLO} for the classification of unipotent classes in $G$. If $q$ is odd, the indecomposable such classes are those labelled $V_\b(2m)$ (a single Jordan block) and $W(m)$ for $m$ odd (two blocks of size $m$); for the latter class, $g = \pm W(m)$ lies in a subgroup $GL_m(q)$ of $G$, and the conclusion follows from Theorem \ref{glgu}. If $q$ is even, the indecomposable unipotent classes are those labelled $V_\b(2m)$, $W(m)$ and $W_\a(m)$ ($m$ odd); elements in the classes $W(m)$ and $W_\a(m)$ lie in subgroups $GL_m(q)$ and $GU_m(q)$ respectively, so again the result follows from Theorem \ref{glgu} in these cases. 

In summary, in case (i) we may take it that $g = \pm V_\b(2m)$, a single Jordan block. The centralizer order is $|C_G(g)| = 2q^m$ (see \cite[Chap. 7]{LS}).

In case (ii), $u$ is a single Jordan block $J_r \in C_G(s)= GL^{\pm}_r(q^d)$ with $rd=m$. The maximal centralizer order occurs when $C_G(s) = GU_m(q)$, and is $(q+1)q^{m-1}$. We conclude that in both cases (i) and (ii),
\begin{equation}\label{spbds}
|C_G(g)| \le 2q^m \hbox{ and } \dim C_V(g) \le 1.
\end{equation}

At this point we use character theory, as in the previous section. Again, we need to show that the sum in (\ref{ndash}) is positive. 

First assume that $q$ is odd. Here we use the following ``gap" result taken from \cite[5.2]{TZ}: $G = Sp_{2m}(q)$ has a collection $\W$ of four irreducible characters of degree $\frac{1}{2}(q^m\pm 1)$ (called Weil characters), such that for all $1 \ne \c \in \hbox{Irr}(G) \setminus \W$, we have $\c(1) \ge \frac{(q^m-1)(q^m-q)}{2(q+1)}$. We break the sum in (\ref{ndash}) into separate parts:
\begin{equation}\label{ndash2}
N' = 1 + \sum_{\c \in\W} \frac{\c(g)^6\c(g^{-1})^6}{\c(1)^{10}} + \sum_{1\ne \c \in Irr(G)\setminus \W} \frac{\c(g)^6\c(g^{-1})^6}{\c(1)^{10}} = 1 + \Sigma_1 + \Sigma_2.
\end{equation}
To bound $\Sigma_1$, we use the following fact (see for example \cite[p.79]{LSh}): for $\c \in \W$, there are complex numbers $a,b$ of modulus 1 such that for all $x \in G$,
\[
\c(x)  = \frac{1}{2}\left(a|C_V(x)|^{1/2} + b|C_V(-x)|^{1/2}\right),
\]
and hence $|\c(g)| \le \frac{1}{2}q^{1/2}$. It follows that 
\[
|\Sigma_1| \le  \frac{4(q^{1/2}/2)^{12}}{((q^m-1)(q^m-q)/2(q+1))^{10}}.
\]
Now we bound $\Sigma_2$. Here we use the bound $k(G) \le 10.8q^m$ from \cite[Prop. 3.12]{FG}. Noting that $|\c(g)| \le |C_G(g)|^{1/2} \le (2q^m)^{1/2}$, we see that
\[
|\Sigma_2| \le  \frac{10.8q^m\cdot (2q^m)^6}{((q^m-1)(q^m-q)/2(q+1))^{10}}.
\]
We now check that for $m\ge 2$ and $q$ odd, with $(m,q) \ne (2,3),$ which we checked using GAP, we have 
$|\Sigma_1| + |\Sigma_2| < 1$. Hence the sum $N'$ in (\ref{ndash2}) is positive, completing the proof of Theorem \ref{sp} when $q$ is odd.

Now assume $q$ is even. Here we use \cite[Thm. 5.5]{TZ}, which states that provided $(m,q) \ne (2,2)$, we have 
 $\c(1) \ge \frac{(q^m-1)(q^m-q)}{2(q+1)}$ for all $\c \in \hbox{Irr}(G) \setminus 1$. Also $k(G) \le 15.2q^m$ by \cite[Thm. 3.13]{FG}, and hence 
\[
 \sum_{\c \in Irr(G)\setminus 1} \frac{|\c(g)^6\c(g^{-1})^6|}{\c(1)^{10}} \le \frac{15.2q^m\cdot (2q^m)^6}{((q^m-1)(q^m-q)/2(q+1))^{10}}.
\]
Once again we check that this is less than 1 for $m\ge 2$ and $q$ even, with $(m,q) \ne (2,2)$ or $(3,2)$ for which cases the result follows by calculations in GAP. %({\bf Please check.}) 

This completes the proof of Theorem \ref{sp}.

\section{Proof of Theorem \ref{othm}} \label{secothm}

This is similar to the previous section, but requires a bit more care in some places. For the induction argument we shall need the following lemma.

\begin{lem}\label{orbase} Let $G = O_4^\e(q)$ and $S = \O_4^\e(q)$, with $q> 3$. Let $g \in G\setminus Z(G)$, and $C = g^S$. Then there exists a non-central involution $t \in S$, and elements $x_i \in C$, $y_i \in C^{-1}$ for $1\le i \le 6$, such that $\prod_{i=1}^6 x_i y_i = t$.
\end{lem}

\pf We have $O_4^\e(q) \cong  \O_4^\e(q).2^a$, where $a=(2,q-1)$, and $\O_4^-(q) \cong PSL_2(q^2)$,   $\O_4^+(q) \cong SL_2(q) \circ SL_2(q)$. Choose $x \in S$ such that $h := [g,x] \in S\setminus Z(S)$. It is sufficient to find 6 conjugates of $h$ with product equal to a non-central involution in $S$. This can be done, as was shown in the proof of Lemma \ref{base}. $\;\;\;\Box$
%({\bf is this right? maybe you can incorporate this lemma with Lemma \ref{base}??}). 

\vspace{4mm}
We now prove Theorem \ref{othm}. Let $G = O_n^\e (q)$ and $S = \O_n^\e(q)$, with $n\ge 5$, and $q$ odd if $n$ is odd. Let $g \in G\setminus Z(G)$, and $C  = g^S$. We assume that 
\begin{equation}\label{oexc}
(n,q) \ne (7,3),\,(8,2),\,(8,3).
\end{equation}
These cases are easily handled by computation. %({\bf Please do it.})

As usual we proceed by induction on $n$. The cases $n=5$ and $n=6$ are covered by Theorems \ref{glgu} and \ref{sp} in view of the isomorphisms $\O_5(q) \cong PSp_4(q)$ and $P\O_6^\e(q) \cong PSL_4^\e(q)$. So assume $n\ge 7$.

\vspace{2mm}
\no {\bf (1) } Consider first the case where $g$ is decomposable, i.e. $g$ lies in a proper subgroup $O_{n_1}(q) \times O_{n_2}(q)$ of $G$ preserving an orthogonal decomposition $V = V_1 \oplus V_2$. Say $g = g_1\oplus g_2$, where $g_i \in O(V_i)$.
We need to be a little careful in applying induction in this case. Let $n_1\ge n_2$, so $n_1 \ge 4$. 

Suppose $g_1 \not \in Z(O(V_1))$. If $n_1 \ge 5$, we can apply induction to find 12 conjugates of $g_1^{\pm 1}$ with product a non-central involution in 
$\O(V_1)$, giving the conclusion in the usual way. If $n_1=4$, then $q>3$ by the exclusions in (\ref{oexc}), and so the conclusion follows in the same way, using Lemma \ref{orbase}.

Now suppose $g_1 \in Z(O(V_1))$, so $g_1 = \d I_{n_1}$ with $\d \in \{1,-1\}$. Take any proper orthogonal decomposition $V_1 = W_1 \oplus W_2$ with $\dim W_1 = 2$. Now replacing $g_1$ by $g_1' = \d I_{W_2} \oplus g_2$, and 
$g_2$ by $\d I_{W_1}$, the argument of the previous paragraph goes through.

\vspace{2mm}
\no {\bf (2) } Now suppose $g$ is indecomposable, and let $g = su$ be the Jordan decomposition. As in the previous section, we have 
\[
\begin{array}{ll}
C_G(s) & = O(V_1) \times O(V_{-1}) \times \prod_{i=1}^k GL_{r_i}(q^{d_i}) \times  \prod_{i=1}^l GU_{s_i}(q^{e_i}) \\
          & \le O(V_1) \times O(V_{-1}) \times \prod_{i=1}^k O_{2r_id_i}(q) \times  \prod_{i=1}^l O_{2s_ie_i}(q).
\end{array}
\]
Since $g$ is indecomposable, there is only one factor in $C_G(s)$, and so one of the following holds:
\begin{itemize}
\item[(i)] $s = \pm I$, $C_G(s) = G$, and $u$ is an indecomposable unipotent element of $G$;
\item[(ii)] $C_G(s) = GL_r(q^d)$ or $GU_r(q^d)$ with $2rd=n$, and $u$ is a single Jordan block $J_r \in C_G(s)$.
\end{itemize}
In case (i), we see as in the previous section, using \cite[Prop. 2.4, Thm. 3.1]{GLO}, that we may take $g$ to be a single Jordan block $\pm V(2m+1)$ ($n=2m+1$, $q$ odd) or $V_\b(2m)$ ($n=2m$, $q$ even). In both cases, $|C_G(g)| \le 2q^m$. In case (ii), as in the previous section we have $|C_G(g)| \le (q+1)q^{m-1}$, where $n=2m$. We conclude that for both (i) and (ii), 
\begin{equation}\label{orpbds}
|C_G(g)| \le 2q^m, \hbox{ where } n=2m \hbox{ or }2m+1.
\end{equation}

We have $g^S = g^{G^*}$ for some $G^*$ such that $S\le G^* \le G$ and $g \in G^*$. In order to apply the usual character-theoretic argument, we need lower bounds on the irreducible character degrees of $G^*$. These are provided by \cite[Thms. 6.1,7.6]{TZ}: for any non-linear irreducible character $\c$ of $G^*$, we have $\c(1) \ge d_G$, where $d_G$ is defined as follows:
\begin{itemize}
\item[(i)] for $G =O_{2m+1}(q)$ with $m\ge 3$, $q$ odd, 
\[
d_G = \left\{ \begin{array}{l} (q^{2m}-1)/(q^2-1), \hbox{ if }q\ge 5 \\ (q^m-1)(q^m-q)/2(q+1), \hbox{ if }q=3 \end{array} \right.
\]
\item[(ii)] for $G =O_{2m}^\e(q)$ with $m\ge 4$ and $(m,q,\e) \ne (4,2,+)$, 
\[
d_G = \left\{ \begin{array}{l} (q^m-\e)(q^{m-1}+\e q)/(q^2-1), \hbox{ if } (q,\e) \ne (2,+),\,(3,+) \\
(q^m-1)(q^{m-1}-1)/(q^2-1), \hbox{ if }(q,\e) = (2,+) \hbox{ or }(3,+) \end{array} \right.
\]
\end{itemize}
Also $k(G^*) \le 15q^m$ by \cite[Thms. 3.14, 3.17, 3.21]{FG}. Using this, together with the bound $\c(g)| \le |C_G(g)|^{1/2} \le (2q^m)^{1/2}$, it follows that 
\[
\sum_{\c \in Irr(G^*), \c(1)>1} \frac{\c(g)^6\c(g^{-1})^6|}{\c(1)^{10}} \le \frac{15q^m\cdot (2q^m)^6}{(d_G)^{10}}.
\]
We check that this is less than 1 unless $(n,q)$ as in the exclusions (\ref{oexc}), and the conclusion of the theorem follows in the usual way.

This completes the proof of Theorem \ref{othm}.

\section{A constructive algorithm}\label{alg}

In this section we give an alternative proof to the case $\epsilon=+1$ of Theorem~\ref{glgu}. Although it requires a larger number of conjugates of $g^{\pm 1}$, the conjugating elements are explicitly determined up to possibly passing to generalized Jordan form, for which algorithms are well-known (see for instance \cite{Steel}). One might be able to apply the method to other groups of Lie type, although with considerably more work.

We fix some notation: $I_{d}$ is the $d\times d$ identity matrix, $E_{ij}$ is the matrix having $ij$-entry $1$ and $0$ elsewhere, and $A\oplus B$ denotes the diagonal join of the matrices $A,B$.

\begin{thm}\label{sl}
Let $G=GL_n(q)$ and $S=SL_n(q)$, with $n\geq 2$ and $(n,q) \neq (2,2),(2,3)$. Let $g \in G\setminus Z(G)$, and $C  = g^S$. Then, for some $k\le 48$, there exists $t\in(CC^{-1})^{k}\setminus Z(G)$ such that $t^{2}\in\{I_{n},-I_{n}\}$.
\end{thm}

%In particular $t$ is an involution in $\mathrm{PSL}_n(q)$ and either $t$ or $t^{2}$ is an involution in $\mathrm{SL}_n(q)$. 

\pf Note first that there exists $h \in S$ such that $[g,h] \in S \setminus Z(S)$. Therefore it is enough to prove the conclusion of the theorem with $g \in S$ and $k\le 24$.

%It is enough to prove the theorem with $g\in S\setminus Z(S)$. In fact every nontrivial class has size $\geq q^{n-1}>(n,q-1)=|Z(S)|$, so $g(g^{-1})^{h}\in S\setminus Z(S)$ for some $h\in S$. Now that we are inside $S$, it is also enough to prove the theorem with $t\in(\{e\}\cup C\cup C^{-1})^{k}$. In fact, for the same reasoning on the size, every nontrivial class $C'$ is such that $C'(C')^{-1}$ contains a nontrivial class $C$ (and the identity $e$, trivially), and if it contains $C$ then it contains $C^{-1}$ as well.

Assume first that $n=2$. For $x \in \F_q$, define $h(x)=I_{2}+xE_{1,2}$. For every $g=(g_{ij})\in SL_{2}(q)$, if $g_{21}\neq 0$ then $g'=h(-g_{11}g_{21}^{-1})gh(-g_{11}g_{21}^{-1})^{-1}$ is such that $g'_{11}=0$. Thus, every class contains an element of at least one of the two forms
\begin{align}\label{eq:sl2start}
g & =\begin{pmatrix}
a & b \\
0 & a^{-1}
\end{pmatrix}, &
g & =\begin{pmatrix}
0 & -a^{-1} \\
a & b
\end{pmatrix}.
\end{align}

%We reduce first to the case $g=h(x)$ with $x\neq 0$ (and thus $g\notin Z(G)$), except for $g=\pm I_{2}$ for which there is nothing to prove. 
Suppose first that $g=h(x)$ with $x\neq 0$. If $q$ is even, this element is an involution already. If $q$ is odd, then there exists $\alpha\in\{-1,2,-2\}$ that is a nonzero square in $\mathbb{F}_{q}$; let $\beta=\frac{2}{\alpha} \in \Z$, and let $\gamma\in\mathbb{F}_{q}^{*}$ satisfy $\gamma^{2}x^{2}=\alpha$. Define
\begin{align*}
s(x) & :=\begin{pmatrix}
-\gamma x & x-\gamma^{-1} \\
\gamma & -1
\end{pmatrix}, & t(x) & :=h(x)s(x)h(x)^{\beta}s(x)^{-1}=\begin{pmatrix}
1 & x \\
-2x^{-1} & -1
\end{pmatrix}.
\end{align*}
Then $t(x)^{2}=-I_{2}$, proving the result in this case (with $k=3$).

Now suppose that $g$ is as in the first case of \eqref{eq:sl2start}, with $a\neq 1$. Define 
\begin{align*}
h_{1} & :=h(1), &
h_{2} & :=h(a^{-1}b), &
h_{3} & :=\begin{pmatrix}
c & (c^{2}-1)ba^{-1}c^{-1} \\
0 & c^{-1}
\end{pmatrix},
\end{align*}
where $c$ is to be chosen below. 
If $a\neq\pm 1$, then $h_{1}gh_{1}^{-1}g^{-1}=h(1-a^{2})$ with $1-a^{2}\neq 0$. And if $a=-1$ then $b\neq 0$, and $g^{2}=h(-2b)$ with $-2b\neq 0$ because $q$ must be odd (or else we would have $a=1$). Hence the result follows from the previous paragraph.

Finally, suppose that $g$ is as in the second case of \eqref{eq:sl2start}. If $b=0$, then $g^{2}=-I_{2}$. If $b\neq 0$ and $q$ is odd, then $(h_{2}gh_{2}^{-1}g)^{2}=h(4a^{-1}b)$ with $4a^{-1}b\neq 0$. If $b\neq 0$ and $q$ is even, take any $c\neq 0$ such that $c^{2}\neq 1$ (note that $c$ exists because $(n,q)\neq(2,2)$ by hypothesis). Then $h_{3}gh_{3}^{-1}g^{-1}$ is as in the first case of \eqref{eq:sl2start} with $c^{2},c^{-2}$ on the diagonal; also $c^{2}\neq \pm 1$, so this case is covered by the previous paragraph, and gives the conclusion with $k=12$. This completes the argument for $n=2$.

Now assume that $n>2$. The cases $(n,q)\in\{(3,2),(3,4),(4,2),(4,3)\}$ can be treated by direct calculation, so exclude these from consideration. Following the proof of Theorem \ref{glgu}, we argue by induction on $n$ and descend by decomposing $g$ or by going to the pairs $(n,q)$ above, unless we are in one of the following situations, corresponding respectively to $m=1,2,n$:
\begin{enumerate}
\item\label{m1} $g=C(f)$, the companion matrix of an irreducible polynomial $f(x)=x^{d}+c_{d-1}x^{d-1}+\ldots+c_{1}x+c_{0}$ of degree $d=n$ (in particular $c_{0}\neq 0$);
\item\label{m2} $g=\begin{pmatrix}C(f)&C(f)\\0&C(f)\end{pmatrix}$, where $f(x)$ has degree $d=\frac{1}{2}n$;
\item\label{mn} $g=\lambda I_{n}+\lambda\sum_{i=1}^{n-1}E_{i,i+1}$ for some $\lambda\in\mathbb{F}_{q}^{*}$.
\end{enumerate}

Start with case \eqref{m1}, and let us do all the calculations as an example. Take
\begin{equation*}
s:=I_{n-2}\oplus \begin{pmatrix}1&-1\\0&1\end{pmatrix},
\end{equation*}
and note that $sg^{-1}s^{-1}=g^{-1}+c_{0}^{-1}E_{n-1,1}+E_{n-2,n}$. Then
\begin{align*}
sg^{-1}s^{-1}g & =(g^{-1}+c_{0}^{-1}E_{n-1,1}+E_{n-2,n})g \\
 & =I_{n}+E_{n-2,n-1}-c_{n-1}E_{n-2,n}-E_{n-1,1} \\
 & =I_{n-3}\oplus\begin{pmatrix}
1 & 1 & -c_{n-1} \\
0 & 1 & -1 \\
0 & 0 & 1
\end{pmatrix} \\
 & =:I_{n-3}\oplus r.
\end{align*}
Then for
\begin{align*}
r(y) & :=\begin{pmatrix}
0 & -1 & y \\
1 & 0 & 0 \\
0 & 0 & 1
\end{pmatrix}, & s(y) & :=I_{n-3}\oplus r(y),
\end{align*}
we get
\begin{equation*}
r(-1)r^{-1}r(-1)^{-1}\cdot r(0)rr(0)^{-1}=\begin{pmatrix}
1 & 0 & 0 \\
0 & 1 & 1 \\
0 & 0 & 1
\end{pmatrix},
\end{equation*}
and putting everything together we obtain in the end
\begin{align*}
\ & s(-1)g^{-1}s(-1)^{-1}\cdot s(0)gs(0)^{-1}\cdot s(1)g^{-1}s(1)^{-1}\cdot s(0)gs(0)^{-1} \\
= \ & s(-1)\cdot(sg^{-1}s^{-1}g)^{-1}\cdot s(-1)^{-1}\cdot s(0)\cdot(sg^{-1}s^{-1}g)\cdot s(0)^{-1} \\
= \ & I_{n-3}\oplus(r(-1)r^{-1}r(-1)^{-1}\cdot r(0)rr(0)^{-1}) \\
= \ & I_{n-2}\oplus\begin{pmatrix}1&1\\0&1\end{pmatrix},
\end{align*}
reducing the problem to the case $n=2$. Note that in the case $n=2$ we obtained the conclusion for $h(1)$ with $k=3$, so it follows here for $g$ with $k=12$.

In case \eqref{m2}, $n$ must be even. Assume $n\geq 8$, take $s(y)$ as in the previous case and
\begin{align*}
t_{1}:= \ & I_{n-2}\oplus\begin{pmatrix}1&1\\0&1\end{pmatrix}, \\
t_{2}:= \ & E_{1,n}+(-1)^{n/2+1}E_{1,n/2+1}+\sum_{i=2}^{n/2-2}E_{i,n/2+i}+\sum_{i=1}^{n/2}E_{n/2-2+i,i} \\
\ & -E_{n-2,n-2}-E_{n-1,n-1}+(1+c_{n/2-1})E_{n-1,n}+E_{n,n}.
\end{align*}
Then
\begin{align*}
h & :=t_{2}^{-1}g^{-1}t_{2}\cdot(t_{2}^{-1}t_{1})g(t_{2}^{-1}t_{1})^{-1}, & s(0)hs(0)^{-1}\cdot s(-1)h^{-1}s(-1)^{-1} & =t_{1},
\end{align*}
reducing again to $n=2$; as above, this gives the conclusion with $k=12$. For $n=6$, use the same procedure replacing
\begin{equation*}
t_{2}:=E_{1,4}+E_{2,3}+E_{2,6}+E_{3,1}+E_{4,2}+E_{6,6}-E_{4,4}-E_{5,5}+(1+c_{2})E_{5,6}.
\end{equation*}
For $n=4$, take instead
\begin{equation*}
t(y):=E_{2,1}+E_{3,3}+E_{4,4}-E_{1,4}-E_{2,1}+yE_{3,4}
\end{equation*}
and
\begin{equation*}
t(0)gt(0)^{-1}\cdot t(-1)g^{-1}t(-1)^{-1}=I_{2}\oplus\begin{pmatrix}1&1\\c_{0}^{-1}&1+c_{0}^{-1}\end{pmatrix},
\end{equation*}
reducing again to $n=2$; in this $n=2$ case we only obtained the conclusion with $k=12$, so it follows for $g$ with $k=24$.

In case \eqref{mn}, define
\begin{equation*}
v(y):=E_{n-1,1}+E_{n-2,2}+\sum_{i=3}^{n}E_{n+1-i,i}+\sum_{i=3}^{n}(-1)^{i+1}E_{n,i}+yE_{n-1,2}.
\end{equation*}
Then
\begin{equation*}
v(1)gv(1)^{-1}\cdot v(0)g^{-1}v(0)^{-1}=I_{2}\oplus\begin{pmatrix}1&1\\0&1\end{pmatrix},
\end{equation*}
reducing again to $n=2$.  $\;\;\;\;\Box$

\section{Proof of Theorem \ref{alt}} \label{secalt}

Let $1 \ne g \in S_n$. The result is easily verified for $n=5$, so assume that $n\ge 6$. 
%We give a constructivese by case proof for this depending on the disjoint cycle decomposition. We will denote the other disjoint cycles of $g$ by $g'.$ \par 
We can write $g = g_0g'$, where $g_0$ is as in the table below, and the cycles of $g'$ are disjoint from those of $g_0$.
In the table, for each possible $g_0$ we provide an element $h \in A_n$ such that $[g,h]$ has order 2, as required. $\;\;\;\Box$

 \[
 \begin{array}{c|c|c}
 g_0 &h&\text{extra conditions}\\
 \hline\hline
(a_1,a_2,\dots,a_{k-1},a_k) &  (a_2,a_5)(a_3,a_6)& k\geq 6\\\hline\hline
(a_1,a_2,a_3,a_4)&(a_1,a_4)(a_2,a_3) &\\\hline\hline
(a_1,a_2,a_3,a_4,a_5)(a_6) &(a_4,a_6,a_5)&\\\hline
(a_1,a_2,a_3,a_4,a_5)(a_6,a_7) &(a_5,a_7,a_6)&\\\hline
(a_1,a_2,a_3,a_4,a_5)(a_6,a_7,a_8) &(a_4,a_8)(a_5,a_6)&\\\hline
(a_1,a_2,a_3,a_4,a_5)(a_6,a_7,a_8,a_9,a_{10}) &(a_4,a_{10})(a_5,a_6)&\\\hline
 \hline
 (a_1,a_2,a_3)(a_4) &(a_2,a_4,a_3)&\\\hline
  (a_1,a_2,a_3)(a_4,a_5) &(a_1,a_3,a_2,a_5,a_4)&\\\hline
  (a_1,a_2,a_3)(a_4,a_5,a_6) &(a_2,a_6)(a_3,a_4)&\\\hline\hline
  (a_1,a_2)(a_3,a_4) &(a_2,a_4,a_3)&\\\hline
   (a_1,a_2)(a_3)(a_4) &(a_1,a_3)(a_2,a_4)&\\\hline\hline
 \end{array}
 \]

\section{A consequence on orbital diameters} \label{orbdiam}

Let $G$ be a finite group acting transitively on set $\O$. 
The orbitals of $G$ are its orbits on $\O\times \O$, and the diagonal orbital is $\{(\a,\a): \a \in \O\}$.
For a non-diagonal orbital $\G$, we define the corresponding orbital graph to be the
undirected graph with vertex set $\O$ and edge set $\{ \{\a,\b\} : (\a,\b) \in \G\}$. 
By \cite[Thm. 3.2A]{DM}, the orbital graphs are all connected if and only if $G$ acts primitively on $\O$,
in which case the {\it orbital diameter} of $G$ is defined to be the supremum of the
diameters of its orbital graphs (see \cite{LMT}). Let us denote this by $orbdiam(G)$.

In \cite{kamilla}, a study is made of the orbital diameters of primitive groups of simple diagonal type.
We are able to use our results to improve one of the theorems in that paper, namely \cite[Thm. 6.1]{kamilla}, as 
follows. Let $T$ be a non-abelian simple group, $k\ge 2$ an integer, and let $S(T,k):= T^k.S_k$ denote the semidirect product 
in which $S_k$ permutes the coordinates in $T^k$ naturally. Define $D = \{(t,\ldots,t) : t \in T\}$, a diagonal subgroup of $T^k$, 
and let $\O = (T^k:D)$ be the set of right cosets of $D$ in $T^k$. Then $S(T,k)$ acts primitively on $\O$, where 
$T^k$ acts by right multiplication, and $S_k$ by permuting the components of coset representatives, 
and this is a primitive group of simple diagonal type (see \cite{LPS}). 

The following result determines the orbital diameter of $S(T,k)$ up to a multiplicative constant. Recall our definition 
$d_T : = \hbox{max}_{C\neq 1} \,\hbox{diam}(\G_C)$ from Section \ref{intr}.

\begin{prop} \label{improve}
Let $G = S(T,k) \cong T^k.S_k$ in the simple diagonal action defined above. Write $C = d(Inv(T))$. Then 
\[
\frac{1}{2}(k-1)d_T \le orbdiam(G) \le 24C(k-1)d_T.
\]
\end{prop}

\pf The lower bound is given by \cite[Thm. 3.1]{kamilla}. An upper bound of $24(k-1)d_T^{\,2}$ is proved in \cite[Thm. 6.1]{kamilla}.
However, in the proof of that result, it is necessary to construct an involution as a product of length $d_T$ of conjugates of some element, which 
accounts for one of the factors $d_T$ in the upper bound. This factor can therefore be replaced by $C$, giving the required upper bound. $\;\;\;\;\Box$

\vspace{2mm}
Note that $d_{T}$ grows linearly in the rank of $T$ for groups of Lie type, and linearly in the degree $n$ of $T=A_{n}$, whereas $C$ is bounded absolutely in all cases. The previous upper bound \cite[Thm. 6.1]{kamilla} was quadratic in $d_T$.

\section*{Acknowledgements}

Daniele Dona was funded by a Young Researcher Fellowship from the Alfr\'ed R\'enyi Institute of Mathematics. Kamilla Rekvényi was funded by the EPSRC Grant EP/W522673/1.

\end{document}